\documentclass[12pt,twoside]{article}

\begin{document}

\centerline{}

\centerline {\Large{\bf Cone Normed Linear Spaces }}
\centerline{}

\newcommand{\mvec}[1]{\mbox{\bfseries\itshape #1}}

\centerline{\bf {T.K. Samanta, Sanjay Roy and Bivas Dinda }}

\centerline{}
\centerline{Department of Mathematics, Uluberia
College, West Bengal, India.}
\centerline{e-mail: mumpu$_{-}$tapas5@yahoo.co.in}

\centerline{Department of Mathematics,}
\centerline{South Bantra Ramkrishna
Institution, West Bengal, India. }
\centerline{e-mail: sanjaypuremath@gmail.com}

\centerline{Department of Mathematics,}
\centerline{Mahishamuri Ramkrishna
Vidyapith, West Bengal, India. }
\centerline{e-mail: bvsdinda@gmail.com}

\centerline{}

\newtheorem{Theorem}{\quad Theorem}[section]

\newtheorem{definition}[Theorem]{\quad Definition}

\newtheorem{theorem}[Theorem]{\quad Theorem}

\newtheorem{remark}[Theorem]{\quad Remark}

\newtheorem{corollary}[Theorem]{\quad Corollary}

\newtheorem{note}[Theorem]{\quad Note}

\newtheorem{lemma}[Theorem]{\quad Lemma}

\newtheorem{example}[Theorem]{\quad Example}

\begin{abstract}
\textbf{\emph{In this paper, we introduce cone normed linear space, study the cone convergence with respect to cone norm. Finally, we prove the completeness of a finite dimensional cone normed linear space.}}
\end{abstract}
${\bf Keywords:}$  \emph{Cone norm, cone convergence, cone continuous, cone complete, cone compact.}

\section{Introduction}
Let $V$ be a real vector space and $(\|\cdot\|)$ be a norm defined on $V$. Our aim is to study cone normed linear spaces. First we go through the definition of Huang and Zhang \cite{Zhang};\\
Let $E$ be a Banach space and $C$ be a non empty subset of $E$. $C$ is called a cone if and only if \\
$(i)\;\;C$ is closed.\\
$(ii)\;\;a,\,b\in \mathbf{R},\,a,b\geq 0,\;x,\,y\in C\;\;\Rightarrow\,ax+by\in C\\$
$(iii)\;\;x\in C\;or\;-x\in C\,$ for every $\,x\in E.\\\\$
Now we define the order relation $'\leq'$  on $E$. For any given cone $C\subseteq E$, we define linearly ordered set $C$ with respect to $'\leq'$ by $x\leq y$ if and only if $y-x\in C.\;x<y\,$ will indicate $x\leq y$ and $x\neq y$; while $x\ll y$ will stands for $y-x\in int C,\; int C$ denotes the interior of $C.\\$
The cone $C$ is called normal cone if there is a number $K>0$ such that for all $x,y \in E,$ \[0\leq x\leq y \,\Leftrightarrow\,\|x\|\leq K\|y\|.\]
The least positive number satisfying the above inequality is called the normal constant of $C.$

\section{Convergence on Cone normed spaces}
Unless otherwise stated throughout this paper we shall denote $\theta$ as the null element of $E$.

\begin{definition}
Let $V$ be a vector space over the field $\mathbf{R}$. The mapping $\|\cdot\|_c:V\longrightarrow E$ is said to be a cone norm if it satisfies the following conditions:\\
$(i) \;\; \|x\|_c\geq \theta\;\;\forall x\in V,\\
(ii) \;\; \|x\|_c= \theta$ if and only if $ x=\theta_V,\\
(iii)\;\;\|\alpha x\|_c=|\alpha|\|x\|\;\;\forall x\in V,\;\alpha\in\mathbf{R},\\ (iv)\;\;\|x+y\|_c\leq \|x\|_c+\|y\|_c\;\;\;\forall x,y\in V.$
\end{definition}

\begin{definition}
Let $V$ be a vector space over the field $\mathbf{R}$. $\|\cdot\|_c$ is a cone norm on $V$. Then $(V,\|\cdot\|_c)$ is called a cone normed linear space.
\end{definition}

\begin{definition}
A sequence $\{x_n\}_n$ in $V$ is said to{\bf convergence} to a point $x\in V$ if for every $\epsilon\in E$ with $\epsilon\gg\theta$ there is a positive integer $n_0$ such that $\|x_n-x\|_c\ll\epsilon,\;\;\;\forall\,n\geq n_0.$ It will be denoted by $\mathop {\lim }\limits_{n\;\, \to
\,\;\infty } x_n=x$ or $x_n\longrightarrow x$ as $n\longrightarrow\,\infty$
\end{definition}

\begin{definition}
A sequence $\{x_n\}_n$ in $V$ is said to be {\bf Cauchy sequence} if for every $\epsilon\in E$ with $\epsilon\gg\theta$ there is a positive integer $n_0$ such that $\|x_n-x_m\|_c\ll\epsilon,\;\;\;\forall\,m,n\geq n_0..$
\end{definition}

\begin{theorem}
Let $(V,\|\cdot\|_c)$ be a cone normed linear space with normal constant $K$. Let $\{x_n\}_n$ be a sequence in $V$. Then $x_n$ converges to $x$ if and only if $\|x_n-x\|_c\longrightarrow \theta.$
\end{theorem}
{\bf Proof.} Let $\{x_n\}_n$ converges to $x$. For every real $p>0$, choose $\epsilon\in E$ with $\epsilon\gg0$ such that $K\|\epsilon\|<p.$ Then there exist a positive integer $n_0$ such that $\|x_n-x\|_c\ll\epsilon,\;\;\forall \,n\geq n_0.\\$
So, $\|\|x_n-x\|_c\|\leq K\|\epsilon\|<p,\;\;\forall\, n\geq n_0.\\ $
Thus $\|x_n-x\|_c\longrightarrow \theta.\\$
Conversely, suppose $\|x_n-x\|_c\longrightarrow \theta$. For $\epsilon\in E$ with $\epsilon\gg0,$ there is a positive number $K\|\epsilon\|$ such that  $\|\|x_n-x\|_c\|\leq K\|\epsilon\|,\;\;\;\forall\, n\geq n_0(\epsilon)$\\
$\Rightarrow \|x_n-x\|_c\leq\epsilon,\;\;\forall \,n\geq n_0$. Hence $x_n$ converges to $x$.

\begin{definition}
The set $B_c(x,r)=\{y\in C:\,\|y-x\|<r\}\;,r\in \mathbf{R}$ is called an open ball in $C$ with center at $x$ and radius $r$.
\end{definition}

\begin{lemma}\label{l2}
Let $x\in C$. Then for every $z\in\,C$, $B_c(x,r)+\{z\}$ is an open ball, Where $B_c(x,r)=\{y\in C:\,\|y-x\|< r\},\;r\in \mathbf{R}.$
\end{lemma}
{\bf Proof.} $B_c(x,r)+\{z\}= \{z\}+\{y\in C:\,\|y-x\|< r\}$
\[=\{z+y:\,\|y-x\|<
r\}=\{p:\,\|p-z-x\|< r\} \hspace{10 cm}\]
\[=\{p:\,\|p-(z+x)\|< r \}=B_c(z+x,r)\hspace{12 cm}\]
Hence the proof.

\begin{corollary}\label{c1}
If $C$ be a cone normed linear space then $C+intC\subseteq intC.$
\end{corollary}

\begin{lemma}\label{l1}
Let $x,y,\epsilon_1,\epsilon_2\in E$ such that $x\ll\epsilon_1$ and $y\ll\epsilon_2$ then $x+y\ll\epsilon_1+\epsilon_2.$
\end{lemma}
{\bf Proof.} Since $x\ll\epsilon_1\,\Rightarrow\,\epsilon_1-x\in intC$ and $y\ll\epsilon_2\,\Rightarrow\,\epsilon_2-y\in intC\\$
Now, $\epsilon_1+\epsilon_2-(x+y)=(\epsilon_1-x)+(\epsilon_2-y)\in intC+intC\subseteq\,C+intC$.\\
Hence by corollary \ref{c1} we have $\epsilon_1+\epsilon_2-(x+y)\in intC.$\\
Thus $x+y\ll\epsilon_1+\epsilon_2.$

\begin{theorem}\label{t1}
In a cone normed linear space $(V,\|\cdot\|_c)$, if $x_n\longrightarrow x$ and $y_n\longrightarrow y$ then $x_n+y_n\longrightarrow x+y$
\end{theorem}
{\bf Proof.} By lemma \ref{l1} the theorem directly follows.

\begin{theorem}\label{t2}
In a cone normed linear space $(V,\|\cdot\|_c)$, if $x_n\longrightarrow x$ and real $\lambda_n\longrightarrow \lambda$ then $\lambda_nx_n\longrightarrow \lambda x.$
\end{theorem}
{\bf Proof.} Obvious.

\begin{theorem}
In a cone normed linear space $(V,\|\cdot\|_c)$, if $\{x_n\}_n$ and $\{y_n\}_n$ are cauchy sequences then $\{x_n+y_n\}_n$ is a cauchy sequence.
\end{theorem}
{\bf Proof.} By lemma \ref{l1} the theorem directly follows.

\begin{theorem}
In a cone normed linear space $(V,\|\cdot\|_c)$, if $\{x_n\}_n$ and $\{\lambda_n\}_n\in \mathbf{R}$ are cauchy sequences then $\{\lambda_n x_n\}_n$ is a cauchy sequence.
\end{theorem}
{\bf Proof.} Obvious.

\begin{definition}
Let $(U,\|\cdot\|_c)$ and $(V,\|\cdot\|_c)$ be two cone normed linear spaces and $f:U\longrightarrow V$ be a function, then $f$ is said to be {\bf cone continuous} at a point $x_0\in U$ if for any given $\epsilon\in E$ with $\epsilon\gg\theta$ there exists $\delta\in E$ with $\delta\gg\theta$ such that $\|x-x_0\|_c\ll\delta\;\Rightarrow\,\|f(x)-f(x_0)\|_c \ll \epsilon.$
\end{definition}

\begin{lemma}\label{l3}
Let $V$ be a  cone normed linear space and $x,y\in V$ then $\|x\|_c-\|y\|_c\leq \|x-y\|_c$ and $\|y\|_c-\|x\|_c\leq \|x-y\|_c.$
\end{lemma}
{\bf Proof.} $\|x\|_c=\|x-y+y\|_c\leq \|x-y\|_c+\|y\|_c.\\$
$\Rightarrow\, \|x-y\|_c+\|y\|_c-\|x\|_c\in C.\\$
$\Rightarrow\, \|x-y\|_c-(\|x\|_c-\|y\|_c)\in C.\\$
$\Rightarrow\, \|x\|_c-\|y\|_c\leq \|x-y\|_c.\\$
Similarly, $\|y\|_c-\|x\|_c\leq \|x-y\|_c.$

\begin{lemma}
Let $x,y,z\in E$ and $x\leq y\ll z$ then $x\ll z$
\end{lemma}
{\bf Proof.} Since $x\leq y\,\Rightarrow\,y-x\in C\,$ and $y\ll z\,\Rightarrow\,z-y\in intC.\\$
Now, $z-x=z-y+y-x\in C+intC.\;$ Then there exist $c_1\in C$ and $c_2\in intC\,$ such that $z-x=c_1+c_2.\;$ Since $c_1\in intC$ then there exists an open ball $B_c(c_1,r)$ such that $c_1\in B_c(c_1,r)\subseteq C.\;$ Therefore $c_1+c_2\in c_2+B_c(c_1,r).\\$
Then by lemma \ref{l2}, $c_1+c_2\in B_c(c_1+c_2,r).\,$ That is $z-x\in B_c(c_1+c_2,r).\\$ Hence $z-x\in intC.$ That is, $x\ll z.$

\begin{theorem}
The cone norm function $f:(U,\|\cdot\|_c)\rightarrow E $ is cone continuous.
\end{theorem}
{\bf Proof.} Let $f(x)=\|x\|_c\;\forall x\in U\,$ and let $x_0\in U.\,$ Since $\|x\|_c,\|x_0\|_c\in E,\\$ So, either $\|x\|_c\geq \|x_0\|_c$ or $\|x\|_c \leq \|x_0\|_c.\\$
So, either $\|x\|_c- \|x_0\|_c\in C$ or $\|x\|_c -\|x_0\|_c\in C.\\$
So, either $\|\|x\|_c-\|x_0\|_c\|_c=\|x\|_c-\|x_0\|_c$ or $\|\|x\|_c-\|x_0\|_c\|_c=\|x_0\|_c-\|x\|_c.\\$
So, by lemma \ref{l3} we have  $\|\|x\|_c-\|x_0\|_c\|_c\leq \|x-x_0\|_c.\\$
Let us choose $\epsilon\in E$ with $\epsilon\gg \theta.\\$ Then, $\|f(x)-f(x_0)\|_c\ll \epsilon\,\Rightarrow \|\|x\|_c-\|x_0\|_c\|_c\leq \|x-x_0\|_c\ll \epsilon.\\$
That is, $\|x-x_0\|\ll \delta\, \Rightarrow\,\|f(x)-f(x_0)\|_c\ll \epsilon\;$ whenever $\delta\in E$ with $\delta\gg \theta$ and $\delta=\epsilon.\\$
Hence cone norm function is cone continuous.

\begin{definition}
A cone normed linear space $(V,\|\cdot\|_c)$ is said to be cone complete if every cauchy sequence in $V$ converges to a point of $V$.
\end{definition}

\begin{theorem}
Let $(V,\|\cdot\|_c)$ be a cone normed linear space such that every cauchy sequence in $V$ has a convergent subsequence then $V$ is cone complete.
\end{theorem}
{\bf Proof.} Let $\{x_n\}_n$ be a cauchy sequence in $V$ and $\{x_{n_k}\}_k$ be a convergent subsequence of $\{x_n\}_n.\;$ Since $\{x_n\}_n$ is a cauchy sequence then for any $\epsilon_1\in E$ with $\epsilon_1\gg \theta$ there exists a positive integer $n_1$ such that \[\|x_n-x_m\|_c\ll \epsilon_1\;\;\;\forall\,m,n\geq n_1\] Let $\{x_{n_k}\}_k$ converges to $x.$ then for any given $\epsilon_2\in E$ with $\epsilon_2\gg \theta$ there exists a positive integer $n_2$ such that \[\|x_{n_k}-x\|_c\ll \epsilon_2\;\;\;\forall\,n_k\geq n_2.\]  Let $n_0=max\{n_1,n_2\}.\\$
$\|x_n-x\|_c=\|x_n-x_{n_k}+x_{n_k}-x\|_c\leq \|x_n-x_{n_k}\|_c+\|x_{n_k}-x\|_c\ll \epsilon_1+\epsilon_2\;\;\;\forall n,n_k\geq n_0$ (by lemma \ref{l1}).\\
Hence the proof.

\begin{theorem}
Let $(V,\|\cdot\|_c)$ be a cone normed linear space and $C$ be a normal cone with normal constant $K$. Then every subsequence of a convergent sequence is convergent to the same limit.
\end{theorem}
{\bf Proof.} Let $\{x_n\}_n$ be a convergent sequence in $V$ and converges to the point $x\in V.$ Let $\{x_{n_k}\}_k$ be a subsequence of $\{x_n\}_n.\\$
Let $\delta\in\mathbf{R}$ then there exist $\epsilon\in E$ with $\epsilon\gg \theta$ such that $K\|2\epsilon\|<\delta.$

Since $x_n$ converges to $x$, then for this $\epsilon\in E\,$ $\exists$ a positive integer $n_0$ such that \[\|x_n-x\|_c\ll \epsilon\;\;\;\;\forall n\geq n_0\]
i.e., $\|x_{n_k}-x\|_c\ll \epsilon\;\;\;\;\forall n_k\geq n_0.$ Hence $x_{n_k}\rightarrow\,x\;\;\;\forall n_k\geq n_0.\\$
If possible let $\{x_{n_k}\}_k$ converges to $y$ also. Then $\exists$ a positive integer $n_1$ such that \[\|x_{n_k}-y\|_c\ll \epsilon\;\;\;\;\forall n_k\geq n_1.\]
Let $n_2=\max\{n_0,n_1\}.$
Now\[\|x-y\|_c=\|x-x_{n_k}+x_{n_k}-y\|_c\leq \|x-x_{n_k}\|_c+\|x_{n_k}-y\|_c\ll \epsilon+\epsilon=2\epsilon\;,\forall n_k\geq n_2\]
Therefore $\|\|x-y\|_c\|\leq K\|2\epsilon\|<\delta.\;$
This implies that $\|\|x-y\|_c\|=0.\\$
Hence the proof.

\section{Finite dimensional cone normed linear spaces}
\begin{lemma}Let $\{x_{1}, x_{2}, . . . , x_{n}\}$ be a linearly independent subset of a cone normed linear space $(V, \|\cdot\|_{c})$. $C$ be a normal cone with normal constant $K$, then there exist an element $c\in int C$ such that for every set of real scalars $\alpha_{1}, \alpha_{2}, . . . , \alpha_{n}$ we have\\ $\|\alpha_{1}x_{1} + \alpha_{2}x_{2} + . . . + \alpha_{n}x_{n}\|_{c}\geq c(|\alpha_{1}|+ |\alpha_{2}| + . . . + |\alpha_{n}|)\cdots\cdots\cdots (1)$
\end{lemma}
\textbf{Proof:} Let $\alpha=|\alpha_{1}|+ |\alpha_{2}| + . . . + |\alpha_{n}|$. If $\alpha=0$ then each $\alpha_{i}$ is zero and hence $(1)$ is true.\\
So we now assume that $\alpha >0$. Then $(1)$ becomes\\
 $\|\beta_{1}x_{1} + \beta_{2}x_{2} + . . . + \beta_{n}x_{n}\|_{c}\geq c\cdots\cdots\cdots (2)$\\
 Where $\beta_{i}= \frac{\alpha_{i}}{\alpha}$ and $\sum_{i=1}^{n}|\beta_{i}|=1.$\\
 It is sufficient to prove that there exists an element $c \in int C$ such that $(2)$ is true for any set of scalars $\beta_{1}, \beta_{2}, . . . , \beta_{n}$ with $\sum_{i=1}^{n}|\beta_{i}|=1$\\
 If possible let this is not true. Then there exists a sequence $\{y_{m}\}_{m}\in V$ where\\
 $y_{m}=  \beta_{1}^{(m)}x_{1} + \beta_{2}^{(m)}x_{2} + . . . + \beta_{n}^{(m)}x_{n}$,  with $\sum_{i=1}^{n}|\beta_{i}^{(m)}|=1$, m=1, 2, ... \\
 such that $\|y_{m}\|_{c}\rightarrow \theta$ as $m\rightarrow \infty$\\
 Since $\sum_{i=1}^{n}|\beta_{i}^{(m)}|=1$ for m= 1, 2, ... . We have $|\beta_{i}^{(m)}|\leq 1$ for i = 1, 2, ... , n; m= 1, 2, ... \\
 Hence for a fixed $i= 1, 2, ... , n$; the sequence $\{\beta_{i}^{(m)}\}_{m}$ is bounded. Therefore by Bolzano-weierstrass theorem $\{\beta_{1}^{(m)}\}_{m}$ has a subsequence converging to $\beta_{1}$(say), and let $\{y_{1,m}\}_{m}$ denote the corresponding subsequence of $\{y_{m}\}_{m}$. By the same reason the sequence
$\{y_{1,m}\}_{m}$ has a subsequence $\{y_{2,m}\}_{m}$(say), for which the corresponding subsequence of real scalars $\{\beta_{2}^{(m)}\}_{m}$ converges to $\beta_{2}$(say). We continuing this process up to n-th stage. At the n-th stage, we obtain a sequence $\{y_{n,m}\} = \{y_{n,1}, y_{n,2}, ... ...\}$ of $\{y_{m}\}_{m}$ whose terms are of the form $y_{n,m}= \sum_{i=1}^{n}\delta_{i}^{(m)}x_{i}$, with $ \sum_{i=1}^{n}|\delta_{i}^{(m)}|=1$, m=1, 2, ... \\
where $\delta_{i}^{m}\rightarrow \beta_{i}$ as $m\rightarrow\infty,\;i=1,2,\cdots,n.$\\
So as $m\rightarrow\infty, y_{n,m}\rightarrow \sum_{i=1}^{n}\beta_{i}x_{i}=y (say)$  where $\sum_{i=1}^{n}|\beta_{i}|=1$ by theorem \ref{t1} and theorem \ref{t2}.\\
This implies that not all $\beta_{i}$ can be zero. Since $\{x_{1}, x_{2}, ... , x_{n}\}$ is linearly independent. Therefore $y\neq \theta_V.$
Now we show that $y_{n,m}\rightarrow y$ implies $\|y_{n,m}\|_{c}\rightarrow \|y\|_{c}$\\
For every real $\epsilon >0$, choose $c\in E$ with $c\gg \theta$ and $K^{2}\|c\|< \epsilon$. Since $y_{n,m}\rightarrow y$ as $m\rightarrow\infty$,then for this element $c$ we can find a positive integer $n_{0}$ such that $\|y_{n,m}-y\|_{c}\ll c, \forall m\geq n_0.$\\
Therefore $\|\|y_{n,m}-y\|_{c}\|\leq K\|c\|,\cdots\cdots\cdots(3)$ \\
since $x\ll y \Rightarrow y-x \in int C\subseteq C \Rightarrow x\leq y$\\
By lemma \ref{l3} we have $\|y_{n,m}\|_{c}-\|y\|_{c}\leq \|y_{n,m}-y\|_{c}$\\
$\Rightarrow \|\|y_{n,m}\|_{c}-\|y\|_{c}\|\leq K\|\|y_{n,m}-y\|_{c}\|\leq K.K\|c\|\;\;\;[$by $(3)]$\\
\smallskip{\hspace{7cm}} $=K^{2}\|c\|< \epsilon,\;\;\; \forall m\geq n_{0}$\\
Hence $\|y_{n,m}\|_{c}\rightarrow \|y\|_{c}$ as $m\rightarrow\infty$\\
As $\{y_{n,m}\}_{m}$ is a subsequence of $\{y_{m}\}_{m}$ and $\|y_{m}\|_{c}\rightarrow \theta$ as $m\rightarrow \infty$\\
Therefore $\|y_{n,m}\|_{c}\rightarrow \theta$ as $m\rightarrow \infty$ and so $\|y\|_{c}=\theta$ which gives $y=\theta_{V}.$\\
This contradiction proves the lemma.

\begin{theorem}
 Every finite dimensional cone normed linear space with normal constant $K$ is cone complete.
\end{theorem}
\textbf{Proof:}Let $(V,\|\cdot\|_c)$ be a cone normed linear space and $C$ be a normal cone with normal constant $K$. Let $\{x_{n}\}$ be an arbitrary cauchy sequence in $V$. We should show that $\{x_{n}\}$ converges to some element $x\in V$. Suppose that the dimension of $V$ is $m$ and let $\{e_{1}, e_{2}, ... , e_{m}\}$ be a basis of $V$. Then each $\{x_{n}\}$ has a unique representation as $x_{n} = \alpha_{1}^{(n)}.e_{1} + \alpha_{2}^{(n)}.e_{2} + ... + \alpha_{m}^{(n)}.e_{m}$ where $\alpha_{1}^{(n)}, \alpha_{2}^{(n)}, ... , \alpha_{m}^{(n)}\in \mathbf{R}$ and $n=1,2,\cdots$\\
Let $\delta\in \mathbf{R}.$ Then there exist $\epsilon\in E$ with $\epsilon\gg \theta$ such that $\frac{\|\epsilon\|}{K\|c\|}<\delta.$
Since $\{x_{n}\}$ is a cauchy sequence, then for this $\epsilon\in E$ there exist a positive integer $n_{0}$ such that $\|x_{n}-x_{r}\|_{c}\ll \epsilon$ for all $n, r \geq n_{0}$\\
By the above lemma it follows that there exist $c\in int C$ such that \\
$\epsilon\gg\|x_{n}-x_{r}\|_{c} = \|\sum_{i=1}^{m}(\alpha_{i}^{(n)}-\alpha_{i}^{(r)})e_{j}\|_{c}\geq c \sum_{i=1}^{m}|\alpha_{i}^{(n)}-\alpha_{i}^{(r)}|$\\
Therefore $\epsilon\gg c \sum_{i=1}^{m}|\alpha_{i}^{(n)}-\alpha_{i}^{(r)}|$\\
Therefore $\|\epsilon\|\geq K\|c \sum_{i=1}^{m}|\alpha_{i}^{(n)}-\alpha_{i}^{(r)}|\|=K \|c\|. \sum_{i=1}^{m}|\alpha_{i}^{(n)}-\alpha_{i}^{(r)}|$\\
Therefore $\sum_{i=1}^{m}|\alpha_{i}^{(n)}-\alpha_{i}^{(r)}|\leq \frac{\|\epsilon\|}{K\|c\|}$\\
$\Rightarrow |\alpha_{i}^{(n)}-\alpha_{i}^{(r)}|\leq \sum_{i=1}^{m}|\alpha_{i}^{(n)}-\alpha_{i}^{(r)}|\leq \frac{\|\epsilon\|}{K\|c\|}<\delta$\\
Therefore $\{\alpha_{i}^{(n)}\}_{n}$ is a cauchy sequence in $\mathbf{R}$ and therefore converges to a real number $\alpha_{i}$, $i=1, 2, ... , m$\\
We now define the element $x= \alpha_{1}e_{1}+\alpha_{2}e_{2}+ ... +\alpha_{m}e_{m}$, which is clearly an element of $V$. Moreover, since $\alpha_{i}^{n}\rightarrow\alpha_{i}$ as $n\rightarrow\infty$ and $i=1, 2, ... , m$\\
We have $\|x_{n}-x\|_{c}=\|\sum_{i=1}^{m}(\alpha_{i}^{(n)}-\alpha_{i})e_{i}\|_{c}\\
\smallskip\hspace{3.5cm}\leq \sum_{i=1}^{m}|\alpha_{i}^{(n)}-\alpha_{i}|\|e_{i}\|_{c}\\
\smallskip\hspace{3.5cm}\rightarrow \sum_{i=1}^{m}0\cdot\|e_{i}\|_{c}\\
\smallskip\hspace{3.5cm}\rightarrow \theta.\\$
This completes the proof.

\begin{definition}
Let $p\in E$ with $\theta\ll p$ and $b\in V$. Define \[B_p(b)=\{x\in V:\,\|x-b\|_c\ll p\}.\]
\end{definition}

\begin{definition}
Let $V$ be a cone normed linear space and $P$ be a subset of $V$. $P$ is said to be closed if for any sequence $\{x_n\}_n$ in $P$ converges to $x\in P.$
\end{definition}

\begin{definition}
Let $V$ be a cone normed linear space. A subset $Q$ of $V$ is said to be the closure of $P(\subseteq V)$ if for any $x\in Q$ there exist a sequence $\{x_n\}_n$ in $P$ such that $x_n\rightarrow x$ as $n\rightarrow \infty$ with respect to the cone norm.
\end{definition}

\begin{definition}
A subset $P$ of a cone normed linear space $V$ is said to be bounded if $P\subseteq B_p(b)$ for some $b\in V$ and $p\in E$ with $\theta\ll p$.
\end{definition}

\begin{definition}
A subset $P$ of a cone normed linear space $V$ is said to be compact if any sequence $\{x_n\}_n$ in $P$, a subsequence can be selected which is convergent to some point of $P$.
\end{definition}

\begin{theorem}
Let $V$ be a cone normed linear space then every cauchy sequence in $V$ is bounded.
\end{theorem}
{\bf Proof.} Let $\{x_n\}_n$ be a cauchy sequence in $V.$\\
$\Rightarrow$ For every $\epsilon\in E$ with $\epsilon\gg \theta$ there exists $n_0\in\mathbf{N}$ such that \[\|x_n-x_m\|_c\ll \epsilon,\;\;\;\forall m,n \geq n_0.\]
$\Rightarrow$ In particular, for $\epsilon=p\in E$ with $p\gg \theta,\;\exists\,k\in \mathbf{N}$ such that \[\|x_n-x_m\|_c\ll p,\;\;\;\forall m,n \geq k.\]
$\Rightarrow\;\;\;\|x_n-x_k\|_c\ll p,\;\;\;\forall n \geq k.\\$
Let $m=maximal\left\{\|x_1-x_k\|_c,\,\|x_2-x_k\|_c,\cdots\cdots,\|x_{k-1}-x_k\|_c\right\}.\\$
Then $\|x_n-x_k\|_c\ll p+m,\;\;\;\forall n \in \mathbf{N}\;\;\Rightarrow\;x_n\in B_{p+m}(x_k).\\$
Hence $\{x_n\}_n$ is bounded.

\begin{theorem}
In a finite dimensional normal cone normed linear space with normal constant $K$, a subset $M$ of $V$ is compact if and only if $M$ is closed and bounded.
\end{theorem}

\textbf{Proof:} Let $M$ be compact subset of $V$. Then by our formal verification it is easy to see that $M$ is closed.\\
Next we show that $M$ is bounded. If possible let $M$ is not bounded. Let $x_{0}\in V$ be a fixed element. Then there exist a point $x_{1}\in M$ such that
$\|x_{1}-x_{0}\|_{c}\geq \epsilon$ for a chosen $\epsilon \in E$ with $\epsilon\gg \theta.$ By the same reason there exists a point $x_{2}\in M$ such that
$\|x_{2}-x_{0}\|_{c}\geq \|x_{1}-x_{0}\|_{c}+ \epsilon$  \\
Continue this process we obtain a sequence $x_{1}, x_{2}, \cdots$ of the set M such that \\
$\|x_{n}-x_{0}\|_{c}\geq \|x_{1}-x_{0}\|_{c}+ \|x_{2}-x_{0}\|_{c}+ \cdots+ \|x_{n-1}-x_{0}\|_{c} + \epsilon \geq \|x_{m}-x_{0}\|_{c} + \epsilon$ for all $m<n$ by lemma 2.16\\
Therefore $\|x_{n}-x_{0}\|_{c} - \|x_{m}-x_{0}\|_{c} \geq \epsilon.\cdots\cdots\cdots (i)$\\
Now $\|x_{n}-x_{0}\|_{c} \leq \|x_{n}-x_{m}\|_{c} + \|x_{m}-x_{0}\|_{c}$\\
Therefore $\|x_{n}-x_{0}\|_{c} + \|x_{m}-x_{0}\|_{c} \leq \|x_{n}-x_{m}\|_{c}\cdots\cdots\cdots (ii)$\\
Using lemma 2.16 we get from $(i)$ and $(ii)$ $\epsilon \leq \|x_{n}-x_{m}\|_{c}$ for all $m<n.$\\
This shows that neither the sequence nor any subsequence of $\{x_{n}\}$ can converge. This contradiction proves that $M$ is bounded.\\
Conversely, let $M$ is closed and bounded and the dimension of $M$ be $n$. Let $\{e_{1}, e_{2}, ... , e_{n}\}$ be a basis of $M$. Let $\{x_{m}\}$ be a sequence in $M$. Since $M$ is bounded, then there exist $p\in E$ such that $x_{n}\in B_{p}(b)$ for some $b\in V$ for all $n\in \mathbf{N}.$ Therefore $\|x_{n}-b\|_{c}\ll p$ for all $n\in \mathbf{N}$.\\
Now $\|x_{n}\|_{c} = \|x_{n}-b+b\|_{c} \leq \|x_{n}-b\|_{c}+\|b\|_{c} \ll p+\|b\|_{c}.$ [by lemma 2.8]\\
Thus $\|x_{n}\|_{c} \ll p+\|b\|_{c}$  for all $n\in \mathbf{N}.$\\
Let $x_{m} = \alpha_{1}^{m}e_{1}+\alpha_{2}^{m}e_{2}+\cdots+\alpha_{n}^{m}e_{n}.$ Where $\alpha_{j}^{m}\in \mathbf{R},$ the set of all real numbers for $m= 1, 2, ... $ and $j= 1, 2, ... , n$. Therefore by lemma 3.1, there exists an element $c\in intC$ such that\\
$p+\|b\|_{c}\gg \|x_{m}\|_{c} = \|\sum_{j=1}^{n}\alpha_{j}^{m}e_{j}\|_{c}\geq c\sum_{j=1}^{n}|\alpha_{j}^{m}|$ for some $c\in intC$ [ by lemma 3.1]\\
Therefore $\|p+\|b\|_{c}\|\geq K\|c\|\sum_{j=1}^{n}|\alpha_{j}^{m}|$
 $\Rightarrow \sum_{j=1}^{n}|\alpha_{j}^{m}|\leq \frac{\|p+\|b\|_{c}\|}{K\|c\|}$ [ as $c\neq \theta$]\\
 Therefore the sequence of numbers $\{\alpha_{j}^{m}\}, m= 1, 2, ... $ and $j= 1, 2, ... , n$ is bounded. So by Bolzano-Weierstrass theorem there exist a convergent subsequence of $\{\alpha_{j}^{m}\}$. Then by the calculation of the lemma 3.1, there exists a subsequence of $\{x_{m}\}$ that converges. Therefore $M$ is compact. Hence the proof.


\begin{thebibliography}{0}
\bibitem{Zhang} H.L.Guang, Z.Xian \emph{ Cone metric spaces and fixed point theorems of contractive mappings}, J. Math. Anal.Appl. $(\,2007\,)$; 332: 1468 - 1476.
\end{thebibliography}
\end{document}